\documentclass[12pt,a4paper]{article}
\usepackage{russ}
\usepackage{graphicx}
\DeclareGraphicsExtensions{.png}
\usepackage{floatflt}
\usepackage{amsfonts}                              %красивые буквы
\usepackage{amsthm}                                %оформление доказательств
\usepackage{amsmath}
\usepackage{gensymb}
\usepackage{amssymb,amsfonts}
\usepackage[usenames]{color}
\usepackage{wrapfig}
\usepackage{url}
\usepackage{cmap}
\usepackage{graphics}
\usepackage{verbatim}
\usepackage[left=1.9cm,right=2cm,
    top=2cm,bottom=2cm,bindingoffset=0cm]{geometry}

\begin{document}

    \theoremstyle{theorem}
    \newtheorem{theorem}{Теорема}
    \newtheorem*{theorem2}{Основная теорема}
    \theoremstyle{definition}
    \newtheorem{definition}{Определение}
    \newtheorem{remark}{Замечание}
    \theoremstyle{lemma}
    \newtheorem{lemma}{Лемма}
    \newtheorem*{lemma prime}{Лемма 4$'$ [3, лемма 8]}
    \theoremstyle{corollary}
    \newtheorem{corollary}{Следствие}
    \theoremstyle{hypothesis}
    \newtheorem{hypothesis}{Гипотеза}
    \theoremstyle{example}
    \newtheorem{example}{Пример}
    \theoremstyle{problem}
    \newtheorem*{question3}{Вопрос}
    \theoremstyle{task}
    \newtheorem*{task}{Задача}
    \newcounter{image}
    \setcounter{image}{1}
    \newcommand{\Q}{\mathbb{Q}}
    \newcommand{\s}{\sqrt{p}}
    \newcommand{\R}{\mathbb{R}}
    \newcommand{\mr}{\leqslant}
    \newcommand{\br}{\geqslant}
%%%%%%%%%%%%%%%%%%%%%%%%%%%%%%%%%%%%%%%%%%%%%%%%%%%%%%%%%%%%%%%%%%%%%%%%%%%%%%%%%%%%%%%%%%%%%%%%%%%%%%%%%%%%%%%%%%%%%%%%%%%%%%%%
%%%%%%%%%%%%%%%%%%%%%%%%%%%%%%%%%%%%%%%%%%%%%%%%%%%%%%%%%%%%%%%%%%%%%%%%%%%%%%%%%%%%%%%%%%%%%%%%%%%%%%%%%%%%%%%%%%%%%%%%%%%%%%%%
%%%%%%%%%%%%%%%%%%%%%%%%%%%%%%%%%%%%%%%%%%%%%%%%%%%%%%%%%%%%%%%%%%%%%%%%%%%%%%%%%%%%%%%%%%%%%%%%%%%%%%%%%%%%%%%%%%%%%%%%%%%%%%%%
%%%%%%%%%%%%%%%%%%%%%%%%%%%%%%%%%%%%%%%%%%%%%%%%%%%%%%%%%%%%%%%%%%%%%%%%%%%%%%%%%%%%%%%%%%%%%%%%%%%%%%%%%%%%%%%%%%%%%%%%%%%%%%%%
%%%%%%%%%%%%%%%%%%%%%%%%%%%%%%%%%%%%%%%%%%%%%%%%%%%%%%%%%%%%%%%%%%%%%%%%%%%%%%%%%%%%%%%%%%%%%%%%%%%%%%%%%%%%%%%%%%%%%%%%%%%%%%%%
%%%%%%%%%%%%%%%%%%%%%%%%%%%%%%%%%%%%%%%%%%%%%%%%%%%%%%%%%%%%%%%%%%%%%%%%%%%%%%%%%%%%%%%%%%%%%%%%%%%%%%%%%%%%%%%%%%%%%%%%%%%%%%%%
    
\date{08.09.2018}

\author{Новиков Иван Викторович \\ Национальный исследовательский университет Высшая школа экономики,\\ факультет математики, г.~Москва}
%\address{Национальный исследовательский университет Высшая школа экономики, факультет математики, г.~Москва}
%\email{iVthik@gmail.com}

\title{Разрезание на подобные прямоугольники многоугольников, составленных из равных прямоугольников.}
\markboth{И.\,В.~Новиков}{Разрезание многоугольника на подобные прямоугольники}
\date{}

\maketitle

    \begin{abstract}
    Пусть многоугольник составлен из равных прямоугольников. Мы находим все квадратичные иррациональности $r$,  при которых рассматриваемый многоугольник разрезается на подобные прямоугольники с отношением сторон $r$.
    \end{abstract}

    Мы рассматриваем следующую задачу: \textit{какие многоугольники можно разрезать на прямоугольники, подобные данному?}

    \smallskip

    Поставленная задача сложна и вряд ли может быть в разумном смысле решена в общем виде. Однако интересны её частные случаи. Мы покажем, как существование разрезания многоугольника на равные прямоугольники помогает установить, разрезается ли многоугольник на подобные прямоугольники с данным отношением сторон. Под \textit{многоугольником} мы понимаем замкнутую связную ограниченную часть плоскости, ограниченную одной или несколькими непересекающимися и несамопересекающимися замкнутыми ломаными. В частности, многоугольник может не быть односвязным. Обозначим $\Q[\s] = \{a + b\s : a, b \in \Q\}$.

    \begin{theorem} %основная теорема
    \label{Main}
        Пусть $a, b, p \in \mathbb{Q}, p > 0, \sqrt{p} \notin \mathbb{Q}$. Пусть многоугольник, все стороны которого принадлежат $\mathbb{Q}[\sqrt{p}]$, разрезан на равные прямоугольники с отношением сторон $y \in \R$. Тогда многоугольник разрезается на подобные прямоугольники с отношением сторон $a + b\sqrt{p}$ тогда и только тогда, когда прямоугольник с отношением сторон $y$ разрезается на подобные прямоугольники с отношением сторон $a + b\sqrt{p}$.
    \end{theorem}

    \begin{corollary}

        Пусть выполнены предположения теоремы \ref{Main}. Тогда:

        \noindent 1) Если $a - b\sqrt{p} > 0$, то многоугольник можно разрезать на прямоугольники с отношением сторон $a + b\sqrt{p}$ тогда и только тогда, когда $$y \in {\{ e + f\sqrt{p} \: : \: e, f \in \mathbb{Q}, e > 0, \frac{|f|}{e} \leqslant \frac{|b|}{a}\}};$$

        \noindent 2) Если $a - b\sqrt{p} < 0$, то многоугольник можно разрезать на прямоугольники с отношением сторон $a + b\sqrt{p}$ тогда и только тогда, когда $$y \in{\{ e + f\sqrt{p} \: : \: e, f \in \mathbb{Q}, f > 0, \frac{|e|}{f} \leqslant \frac{|a|}{b}\}}.$$

    \end{corollary}

    \begin{wrapfigure}{o}{0.2\linewidth} %рис.1
    \center{\includegraphics[width=1.0\linewidth]{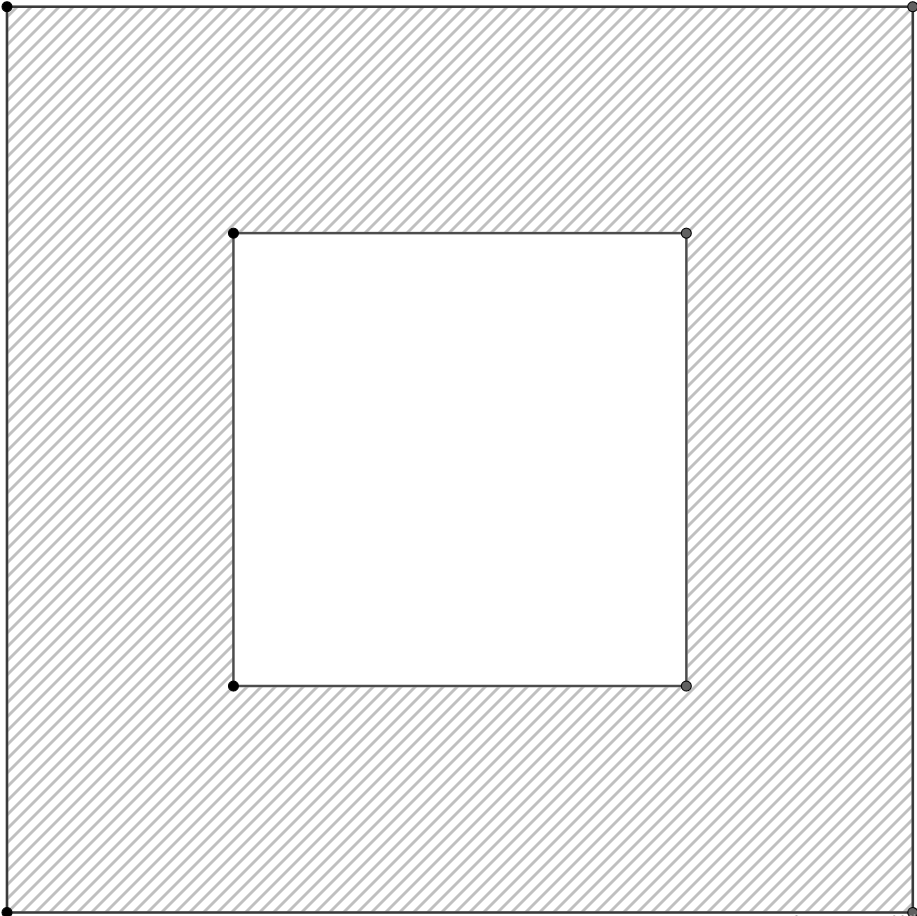} \\ рис. \theimage}
    \end{wrapfigure}

    В качестве ещё одного следствия рассмотрим разрезание "квадрата с дыркой" (рис. \theimage).
    \stepcounter{image}

    \smallskip

    \begin{corollary} %квадрат с дыркой
        Пусть $a, b, p \in \Q, \s \notin \Q$. Рассмотрим многоугольник между двумя гомотетичными концентрическими квадратами со сторонами $u > v$. Тогда:

        \noindent 1) Если $a - b\s > 0$, то данный многоугольник можно разрезать на прямоугольники с отношением сторон $a + b\s$ тогда и только тогда, когда $$\frac{u + v}{u - v} \in{\{ e + f \s \: : \: e, f \in \Q, e > 0, \frac{|f|}{e} \mr \frac{|b|}{a}\}};$$

        \noindent 2) Если $a - b\s < 0$, то данный многоугольник можно разрезать на прямоугольники с отношением сторон $a + b\s$ тогда и только тогда, когда $$\frac{u + v}{u - v} \in{\{ e + f \s \: : \: e, f \in \Q, f > 0, \frac{|e|}{f} \mr \frac{|a|}{b}\}}.$$
    \end{corollary}

    Рассмотрим кратко историю вопроса о разрезании прямоугольника на подобные прямоугольники. Более подробно об этом можно прочитать, например, в \cite['Обзор известных результатов']{3}. Одним из первых нетривиальных результатов был доказанный М. Деном в 1903 г. факт, что квадратами можно замостить прямоугольники только с рациональным отношением сторон \cite{1}. В 2016 г. в \cite{3} Ф. Шаров обобщил теорему Дена, найдя все прямоугольники, которые можно разрезать на прямоугольники с данными отношениями сторон из $\Q[\s]$. Благодаря теореме Шарова мы получим теорему \ref{Main}, а также следствия 1 и 2. В статье \cite{2} Фрайлинг, Лацкович и Ринн нашли алгебраический критерий возможности разрезания прямоугольника на подобные данному прямоугольники. Однако данный критерий не даёт алгоритма проверки возможности данного разрезания.

    %Рассмотрим кратко историю вопроса о разрезании прямоугольника на подобные прямоугольники. Более подробно об этом можно прочитать, например, в \cite[см. 'Обзор известных результатов']{3}. Один из первых нетривиальных результатов, связанных с разрезанием фигуры на подобные прямоугольники, был доказан М. Деном в 1903 Было доказано, что квадратами можно замостить прямоугольники только с рациональным отношением сторон. В 2016 г. Ф. Шаровым было получено обобщение теоремы Дена: были найдены все прямоугольники, которые можно разрезать на прямоугольники с данными отношениями сторон из $\Q[\s]$. Благодаря теореме Шарова из теоремы \ref{Main} выводится следствие 1. В статье \cite{2} Фрайлинг, Лацкович и Ринн нашли алгебраический критерий возможности разрезания прямоугольника на подобные прямоугольники. Однако неизвестно алгоритма проверки данного критерия. %При доказательстве теоремы \ref{Main} будут рассматриваться разрезания многоугольников на равные прямоугольники. Случай разрезания прямоугольников на равные прямоугольники полностью исследован в \cite{10}.

    Зафиксируем число $p \in \Q$ такое, что $p > 0, \s \notin \Q$. Будем называть многоугольник \textit{хорошим}, если его стороны имеют вид $\alpha + \beta\s$ для некоторых $\alpha, \beta \in \Q$. Начиная с этого момента будем предполагать, что стороны всех рассматриваемых многоугольников параллельны координатным осям, т.е. \textit{вертикальны} и \textit{горизонтальны}.

    Мы используем следующие понятия, аналогичные введённым в \cite{13}, ср. с \cite{2}.

    \begin{definition} %базис
        \textit{Базисом} будем называть упорядоченный набор вещественных чисел, линейно независимый над полем $\Q[\s]$. Будем говорить, что число \textit{записывается} в базисе ($e_1, e_2, \ldots, e_k$), если оно равно некоторой линейной комбинации чисел $e_1, e_2, \ldots, e_k$ с коэффициентами из $\Q[\s]$.
    \end{definition}

    \begin{definition}
        Пусть $E = (e_1, e_2, \ldots, e_k)$ — некоторый базис и $z \in \R$. Будем называть \textit{z-площадью в базисе $Е$} (или \textit{площадью Гамеля}) прямоугольника $$(a_1 e_1 + a_2 e_2 + \ldots + a_k e_k) \times (b_1 e_1 + b_2 e_2 + \ldots + b_k e_k),$$
        число $(a_1 + a_2 z) (b_1 + b_2 z)$, где $a_1, a_2, \ldots, a_k, b_1, b_2, \ldots, b_k \in \Q[\s]$.
    \end{definition}

    \begin{lemma}
    \label{y-area}
        Если прямоугольник разрезан на прямоугольники со сторонами, которые записываются в базисе $E$, то $z$-площадь в базисе $E$ разрезаемого прямоугольника равна сумме $z$-площадей в базисе $E$ прямоугольников, на которые он разрезан.
    \end{lemma}

    \begin{wrapfigure}[10]{o}{5.5cm}
    \center{\includegraphics[width=0.85\linewidth]{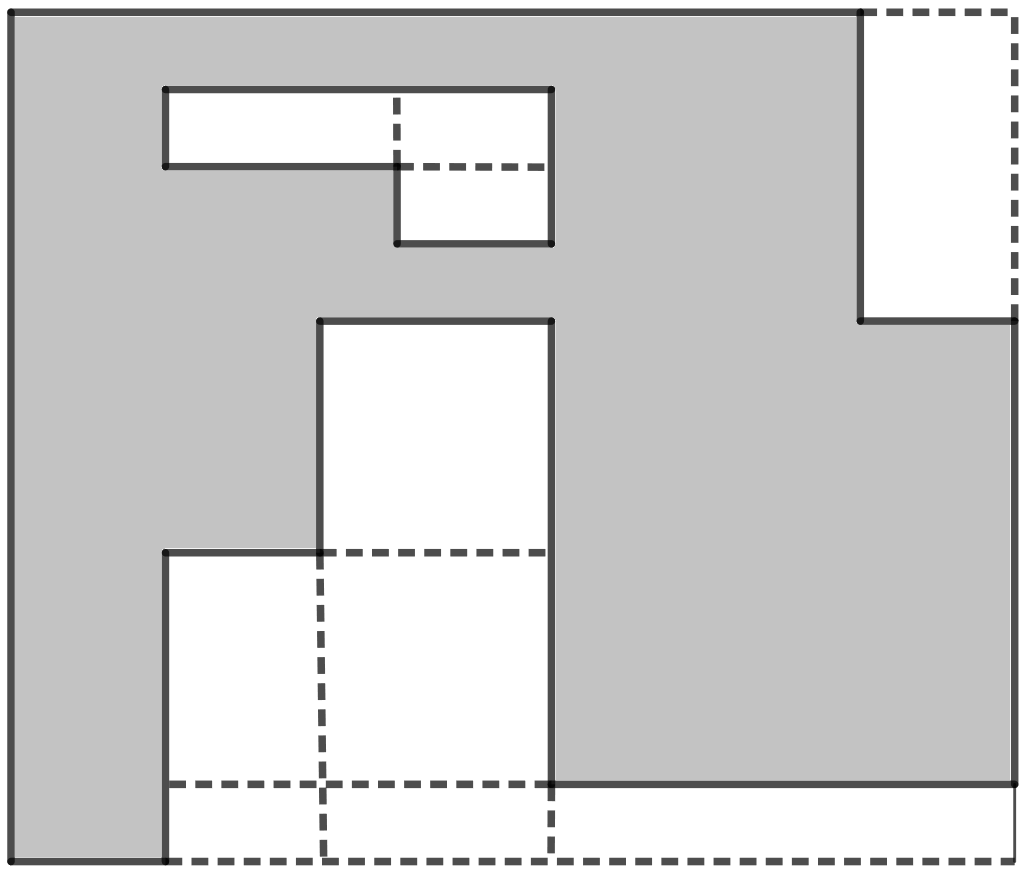} \\ рис. \theimage}
    \end{wrapfigure}

    Доказательство данного свойства $z$-площадей почти дословно повторяет доказательство леммы 8 из \cite{3}.

    \begin{lemma}%[(см. рис. \theimage)]
    \label{Construction}
        (см. рис. \theimage).
        Для любого хорошего многоугольника $M$ с углами, кратными $90 \degree$, найдутся такие хорошие прямоугольники $R, R_1, \ldots, R_n$, что $\text{R}$ разрезается на $M, R_1, \ldots, R_n$.
        \stepcounter{image}
    \end{lemma}

    \begin{proof}
        Искомый прямоугольник $R$ получается как пересечение продолжений самой высокой и низкой горизонтальных сторон многоугольника $M$ с его самой левой и правой вертикальными сторонами. Для получения остальных искомых прямоугольников $R_1, \ldots, R_n$ продолжим каждую сторону многоугольника до пересечения с границей многоугольника $M$ или прямоугольника $R$. Таким образом у нас получится хороший прямоугольник $R$ вместе с разрезанием на хорошие прямоугольники $R_1, \ldots, R_n$ и исходный многоугольник $M$. 
    \end{proof}

    \begin{lemma} %равные прямоугольники
    \label{Equal}
        Если хороший многоугольник разрезан на равные прямоугольники, то эти равные прямоугольники являются хорошими.
    \end{lemma}
    \begin{proof}
        Пусть имеется разрезание хорошего многоугольника на $N$ равных прямоугольников $c \times d$. Пусть, от противного, одно из чисел $c$ или $d$ не принадлежит $\Q[\s]$. Площадь хорошего многоугольника $N c d \in \Q[\s]$, следовательно, если $d \in \Q[\s]$, то и $c \in \Q[\s]$. Значит, $c \notin \Q[\s]$ и $d \notin \Q[\s]$. Рассмотрим базис
        \begin{equation*}
            E =
                \begin{cases}
                    (c, d, 1),&\text{если тройка $(c, d, 1)$ линейно независима над $\Q[\s]$}; \hspace{58pt} (1)\\
                    (1, c),&\text{если тройка $(c, d, 1)$ линейно зависима над $\Q[\s]$}.$ \hfill (2)$
                \end{cases}
        \end{equation*}
        В случае (2) пусть $d = a_1 + a_2 c$ для некоторых $a_1, a_2 \in \Q[\s]$. Тогда $a_2 \neq 0$, так как в противном случае $d$ принадлежало бы множеству $\Q[\s]$.

        По лемме \ref{Construction} дополним многоугольник до хорошего прямоугольника $R$. Получим разрезание хорошего прямоугольника $R$ на хорошие прямоугольники и прямоугольники $c \times d$. Тогда $z$-площадь в базисе $E$ любого хорошего прямоугольника не зависит от $z$. Значит, $z$-площадь прямоугольника $R$ не зависит от $z$. Сумма $z$-площадей прямоугольников, на которые он разрезан, равна $N z$ в случае (1) и $C + N a_1 z + N a_2 z^2$, где $C$ — независящее от $z$ число, в случае (2). По лемме \ref{y-area} получаем, что $z$-площадь прямоугольника $R$ зависит от $z$, так как $a_2 \neq 0$ и $N \neq 0$. Противоречие. 
    \end{proof}

    \begin{lemma}
    \label{Ratio -> sides}
        Если хороший многоугольник разрезан на прямоугольники, отношения сторон которых принадлежат множеству $\Q[\s]$, то эти прямоугольники являются хорошими.
    \end{lemma}
    \begin{proof}
        Пусть, от противного, в разрезании есть прямоугольник $c \times k c$, где $k \in \Q[\s], k > 0, c \notin \Q[\s]$. Дополним набор $(1, c)$ до базиса следующим образом. Рассмотрим множество всех сторон прямоугольников разрезания, и станем по очереди добавлять к нашему набору длины тех сторон, которые оставляют набор линейно независимым над $\Q[\s]$. В итоге получим некоторый базис $E = (1, c, e_3, \ldots, e_n)$. Тогда $z$-площадь произвольного прямоугольника разрезания $$(a_1 + a_2 c + a_3 e_3 + \ldots + a_n e_n) \times b (a_1 + a_2 c + a_3 e_3 + \ldots + a_n e_n),$$
        где $b \in \Q[\s], b > 0$, в базисе $E$ равна $b a_1^2 + 2 b a_1 a_2 z + b a_2^2 z^2$. Коэффициент при $z^2$ неотрицательный для любого прямоугольника разрезания, и положительный для прямоугольника $c \times k c$. Значит, сумма $z$-площадей всех прямоугольников разрезания имеет нетривиальную зависимость от $z$. По лемме \ref{Construction} дополним многоугольник до хорошего прямоугольника. Тогда $z$-площадь полученного прямоугольника, а также добавленных прямоугольников, не зависит от $z$. Значит, по лемме \ref{y-area} получаем противоречие. 
    \end{proof}

    Нам понадобится ещё одно понятие из работы \cite{3}.

    \begin{definition} %площадь прямоугольника
        \cite[определение 5]{3}. Пусть $A,B,C \in \R$,  $\alpha, \beta, \gamma, \delta, p \in \Q$ и $\s \notin \Q$. Тогда \textit{$ABC$-площадь прямоугольника} $(\alpha + \beta\s) \times (\gamma + \delta\s)$ — это число
        $$S := \alpha\gamma A + \alpha\delta B + \beta\gamma B + \beta\delta C.$$
    \end{definition}

    \begin{example}
        Пусть $p=2, A = -1, B = 1, C \in \R$. Тогда $ABC$-площадь прямоугольника $(\sqrt{2} + 1) \times 1$ равна 0, а $ABC$-площадь подобного ему прямоугольника $(\sqrt{2} - 1) \times 1$ равна 2.
    \end{example}

    \begin{lemma}
    \label{ABC-area rectangle}
        \cite[лемма 8]{3}.
        Если прямоугольник разрезан на хорошие прямоугольники, то $ABC$-площадь разрезаемого прямоугольника равна сумме $ABC$-площадей прямоугольников, на которые он разрезан.
    \end{lemma}

%    Приведённое в \cite{3} доказательство по сути ничем не отличается от доказательства леммы \ref{y-area} выше.
    \begin{definition} %площадь многоугольника
        Пусть существует хотя бы одно разрезание данного хорошего многоугольника на хорошие прямоугольники. Тогда \textit{$ABC$-площадь многоугольника} положим равной сумме $ABC$-площадей всех прямоугольников любого такого разрезания.
    \end{definition}

    Докажем корректность определения $ABC$-площади многоугольника.

    \begin{lemma}
    \label{ABC-area polygon}
        Если многоугольник двумя способами разрезан на хорошие прямоугольники, то сумма $ABC$-площадей прямоугольников в обоих способах одна и та же.
    \end{lemma}

    \begin{proof}
         По лемме \ref{Construction} дополним многоугольник до хорошего прямоугольника. По лемме \ref{ABC-area rectangle} для образовавшегося разрезания хорошего прямоугольника $ABC$-площадь равна сумме $ABC$-площадей добавленных хороших прямоугольников и $ABC$-площадей прямоугольников, на которые разрезан многоугольник. Значит, сумма $ABC$-площадей прямоугольников разрезания равна разности $ABC$-площадей прямоугольника и добавленных хороших прямоугольников, а для них $ABC$-площадь определена однозначно и не зависит от разрезания. Значит, и $ABC$-площадь многоугольника не зависит от разрезания. 
    \end{proof}

%    \begin{lemma}
%    \label{Ratio -> sides}
%        Если хороший многоугольник разрезан на конечное число прямоугольников с отношением сторон $a + b\s$, то стороны этих прямоугольников принадлежат множеству $\Q[\s]$.
%    \end{lemma}
%    \begin{proof}
%        По лемме \ref{Construction} дополним многоугольник до хорошего прямоугольника. Получим разрезание хорошего прямоугольника на хорошие прямоугольники. По теореме \ref{Tok} о силах тока и длинах сторон и лемме \ref{Chain} о сопротивлении цепи стороны прямоугольников разрезания выражаются через сторону прямоугольника, полученного дополнением многоугольника прямоугольниками, и отношения сторон хороших прямоугольников с помощью сложения, вычитания, умножения и деления. Тогда, поскольку множество $\Q[\s]$ замкнуто относительно этих операций, то стороны прямоугольников разрезания также принадлежат этому множеству.
%    \end{proof}

        Сформулируем лемму, из которой мы выведем теорему \ref{Main}.

        \begin{lemma}
        \label{Main lemma}
            Пусть $a, b, e, f, p \in \Q, \s \notin \Q, e + f\s > 0, a + b\s > 0$. Предположим, что выполнено одно из двух условий:

            \noindent а) \noindent $a - b\s > 0, \frac{|f|}{e} > \frac{|b|}{a}$;

            \noindent б) \noindent $a - b\s < 0, \frac{|e|}{f} > \frac{|a|}{b}$.

            \noindent Положим в определении $ABC$-площади $A = f, B = -e, C = \frac{2fa^2}{b^2} - pf$. Тогда:

            \noindent 1) $ABC$-площадь прямоугольника со сторонами $e + f\s$ и 1 равна нулю.

            \noindent 2) $ABC$-площадь всех прямоугольников с отношением сторон $a + b\s$ и сторонами из $\Q[\s]$ отличны от нуля и имеют один и тот же знак.
        \end{lemma}

        \noindent \textit{Доказательство теоремы \ref{Main}, считая лемму \ref{Main lemma} доказанной.}
          Докажем часть "только тогда". Рассмотрим два разбиения многоугольника: первое — на равные прямоугольники, второе — на подобные хорошие прямоугольники. Тогда согласно лемме \ref{Equal} и лемме \ref{Ratio -> sides} стороны всех этих прямоугольников принадлежат множеству $\Q [\s]$. Без ограничения общности, стороны равных прямоугольников первого разрезания будем считать равными $e + f\s$ и $1$ (чего можно добиться с помощью подходящей гомотетии). Если числа $a,b,p$ таковы, что оба условия а) и б) леммы \ref{Main lemma} не выполнены, то по части "тогда" в теореме 1 из \cite{3} мы получаем требуемое разрезание прямоугольника с отношением сторон $y$ (прямое построение разрезания см. в \cite[\S 1]{4}). В противном случае по пункту 1) леммы \ref{Main lemma} $ABC$-площадь многоугольника равна нулю (как сумма $ABC$-площадей равных прямоугольников). Но по пункту 2) леммы \ref{Main lemma} сумма $ABC$-площадей прямоугольников с отношениями сторон $a + b\s$ всегда не ноль. Следовательно, $ABC$-площадь многоугольника не равна нулю. Противоречие.

        \smallskip

        Осталось доказать лемму \ref{Main lemma}. Она фактически доказана в \cite[стр 209-211, см. 'Доказательство основной теоремы в случае n = 1.']{3}, но в явном виде там не сформулирована. Приведём здесь её доказательство, по сути заимствованное из \cite{3}.

        \begin{proof}[леммы \ref{Main lemma}]
            1) $ABC$-площадь прямоугольника $(e + f\s) \times 1$ равна

            \begin{center}
                $S = e \cdot 1 \cdot f + e \cdot 0 \cdot (-e) + f \cdot 1 \cdot (-e) + f \cdot 0 \cdot (\frac{2fa^2}{b^2} - pf) = 0$.
            \end{center}

            2) Рассмотрим прямоугольник $(\alpha + \beta\s) \times (\gamma + \delta\s)$, где

            \begin{center}
            $\frac{\gamma + \delta\s}{\alpha + \beta\s} = a + b\s$.
            \end{center}

            Тогда $\gamma + \delta\s = \alpha a + p \beta b + \s(\beta a + \alpha b)$, и значит, $\gamma = \alpha a + p \beta b, \delta = \beta a + \alpha b$. Поэтому $ABC$-площадь прямоугольника записывается в виде
            \begin{align*}
            S &= \alpha \gamma f - \beta \gamma e - \alpha \delta e + \beta \delta \left(\frac{2fa^2}{b^2} - p f\right) = \\
            &= \alpha f (\alpha a + p \beta b) - \beta e (\alpha a + p \beta b) - \alpha e (\beta a + \alpha b) + \beta \left(\frac{2fa^2}{b^2} - p f\right)(\beta a + \alpha b) = \\
            &= {\alpha}^2 (f a - e b) + 2 \alpha \beta \left(\frac{fa^2}{b} - e a\right) + {\beta}^2 \left(\frac{2fa^3}{b^2} - p f a - p e b\right).
            \end{align*}

            Заметим, что при выполнении условия а) мы имеем $a - b\s > 0$ и $a + b\s > 0$, значит $a > 0$. Так как в этом случае $\frac{|f|}{e} > \frac{|b|}{a}$, то $e > 0$. Имеем
            $$\frac{|f|}{e} \neq \frac{|b|}{a} \Leftrightarrow \biggl|\frac{f}{e}\biggr| \neq \biggl|\frac{b}{a}\biggr| \Leftrightarrow |fa| \neq |be| \Rightarrow fa \neq be.$$

            При выполнении условия б) аналогично получаем, что $f a - e b \neq 0$. Значит, $f a - e b \neq 0$ и в случае а), и в случае б). Следовательно, выражение $S / \beta ^2$ — квадратный трёхчлен относительно $\alpha / \beta$. Покажем, что его дискриминант $D$ отрицателен. Тем самым будет доказано, что при фиксированных числах $a, b, e, f$ величина $S$ либо всегда положительна, либо всегда отрицательна.

            Действительно,
            \begin{align*}
            \frac{D}{4} &= \frac{f^2 a^4}{b^2} - \frac{2 f e a^3}{b} + e^2 a^2 - \frac{2 f^2 a^4}{b^2} + p f^2 a^2 + p f a e b + \frac{2 f e a^3}{b} - p f a e b - p e^2 b^2 = \\
            &= \frac{-f^2 a^4}{b^2} + e^2 a^2 + p f^2 a^2 - p e^2 b^2 = (a^2 - p b^2)\left(e^2 - \frac{f^2 a^2}{b^2}\right) < 0,
            \end{align*}
            так как при условии а) из формулировки леммы первая скобка положительна, а вторая отрицательна, а при условии б) — наоборот.

            Таким образом, утверждение 2), а с ним и теорема \ref{Main}, доказано. 
        \end{proof}

    \begin{wrapfigure}[8]{o}{0.26\linewidth}
    \center{\includegraphics[width=1.0\linewidth]{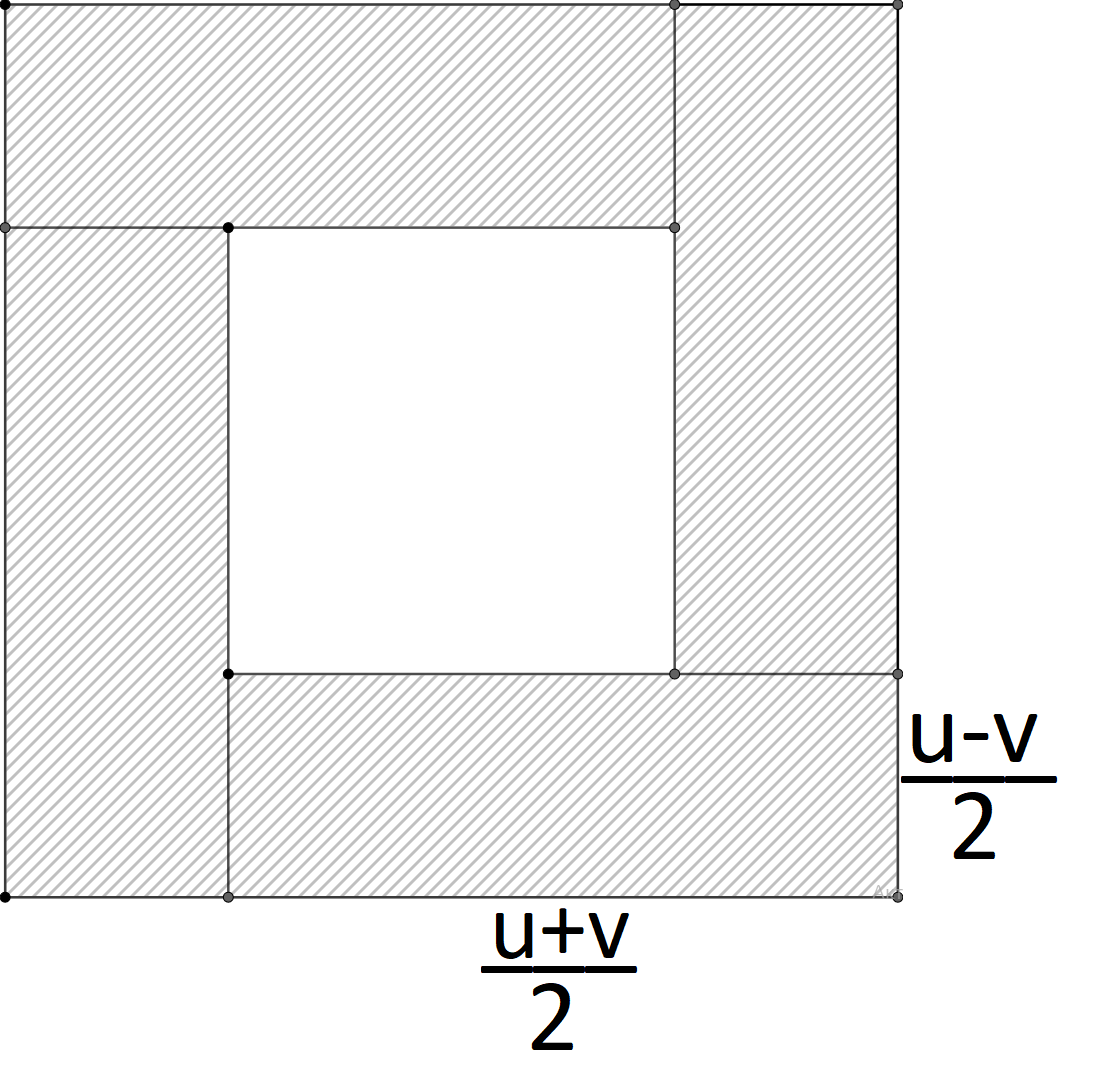} \\ рис. \theimage}
    \end{wrapfigure}
    Докажем сформулированные выше следствия. \stepcounter{image}

    \begin{proof}[Доказательство следствия 1]
        Сразу следует из теоремы \ref{Main} и \cite[Теорема 1]{3}. 
    \end{proof}

    \begin{wrapfigure}[14]{o}{0.52\linewidth}
    \center{\includegraphics[width=0.98\linewidth]{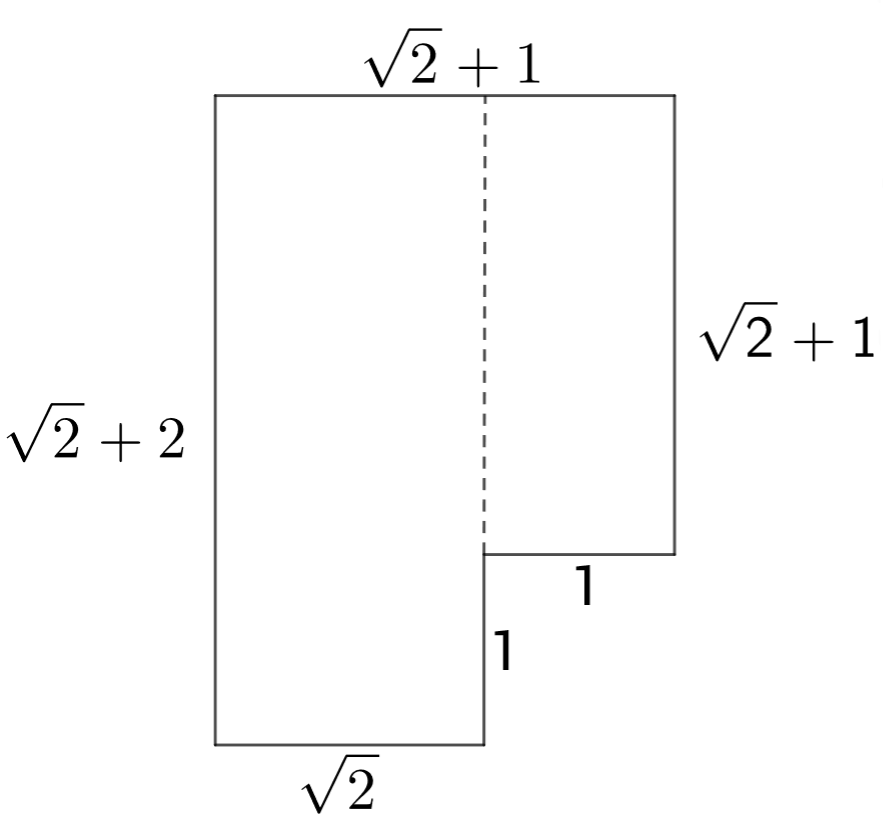} \\ рис. {\theimage}}
    \addtocounter{image}{-1}
    \end{wrapfigure}
    \begin{proof}[Доказательство следствия 2]
        Разрежем квадрат с дыркой на 4 равных прямоугольника $\frac{u + v}{2} \times \frac{u - v}{2}$, как на рис. \theimage. Положим $u = 1$ без ограничения общности.

        Если $\frac{u + v}{u - v} \in \Q[\s]$, квадрат с дыркой является хорошим многоугольником, и по следствию 1 получаем следствие 2.

        Если $\frac{u + v}{u - v} \notin \Q[\s]$, то следствие 1 неприменимо, так как многоугольник не хороший. Докажем, что в этом случае разрезания квадрата с дыркой на прямоугольники с отношением сторон $a + b\s$ не существует. Пусть существует разрезание на прямоугольники с отношением сторон $a + b\s$. Добавляя к разрезанию квадрат $v \times v$, получим разрезание квадрата $u \times u$. По лемме \ref{Ratio -> sides} стороны всех прямоугольников разрезания (включая квадрат $v \times v$) лежат в множестве $\Q[\s]$. Тогда $v \in \Q[\s]$ и $\frac{u + v}{u - v} \in \Q[\s]$. Противоречие. \stepcounter{image} 
    \end{proof}

    %\stepcounter{image}
    Теорема \ref{Main} не даёт универсального метода решения задачи, сформулированной в начале статьи. Например, есть хороший многоугольник, который можно разрезать на подобные хорошие прямоугольники, но нельзя на равные. На рис. \theimage~показан такой многоугольник и его разрезание на подобные прямоугольники. Доказательство отсутствия разрезания на равные прямоугольники опустим.

    В заключение сформулируем в виде задачи связанные с теоремой \ref{Main} открытые вопросы.

    \begin{task}
    Остается ли теорема \ref{Main} справедливой, если в ней заменить "$a + b\s$":

    1) На набор $a_1 + b_1\s, \ldots, a_n + b_n\s \in \Q[\s]$ (которому должны принадлежать отношения сторон прямоугольников)?

    2) На корень $r$ многочлена степени 3 с целыми коэффициентами (а $\Q[\s]$ заменить на $\Q[r]$)?

    3) На произвольное $r \in \R$ (а $\Q[\s]$ заменить на $\R$)?
    \end{task}

    \smallskip

    Автор выражает благодарность М.Б. Скопенкову за постоянное внимание к данной работе.

\end{document}